\documentclass[aps,pre,floats,floatfix,twocolumn,superscriptaddress]{revtex4}
\usepackage{latexsym,amsmath}
\usepackage{amsbsy}
\usepackage[pdftex]{graphicx}
\usepackage{amssymb}
\usepackage{epstopdf}
\usepackage[usenames]{color}

\begin{document}

\title{The complex architecture of primes and natural numbers}

\author{Guillermo Garc\'{i}a-P\'{e}rez}
\author{M. \'{A}ngeles Serrano}
\author{Mari\'{a}n Bogu\~{n}\'{a}}
\affiliation{Departament de F\'{i}sica Fonamental, Universitat de Barcelona \\ Mart\'{i} i Franqu\`{e}s 1, 08028 Barcelona, Spain}

\date{\today}

\begin{abstract}
Natural numbers can be divided in two non-overlapping infinite sets, primes and composites, with composites factorizing into primes. Despite their apparent simplicity, the elucidation of the architecture of natural numbers with primes as building blocks remains elusive. Here, we propose a new approach to decoding the architecture of natural numbers based on complex networks and stochastic processes theory. We introduce a parameter-free non-Markovian dynamical model that naturally generates random primes and their relation with composite numbers with remarkable accuracy. Our model satisfies the prime number theorem as an emerging property and a refined version of Cram\'er's conjecture about the statistics of gaps between consecutive primes that seems closer to reality than the original Cram\'er's version. Regarding composites, the model helps us to derive the prime factors counting function, giving the probability of distinct prime factors for any integer.  Probabilistic models like ours can help to get deeper insights about primes and the complex architecture of natural numbers.

\end{abstract}

\maketitle

\section{Introduction}
Prime numbers have fascinated and puzzled philosophers, mathematicians, physicists and computer scientists alike for the last two and a half thousand years. A prime is a natural number that has no divisors other than $1$ and itself; every natural number greater than $1$ that is not a prime is called a composite. Despite the apparent simplicity of these definitions, the hidden structure in the sequence of primes and their relation with the set of natural numbers are not yet completely understood~\footnote{Leonhard Euler about the structure of primes: \textit{``Mathematicians have tried in vain to this day to discover some order in the sequence of prime numbers, and we have reason to believe that it is a mystery into which the human mind will never penetrate.''} ~---~ Leonhard Euler, 1751}. There is no practical closed formula that sets apart all of the prime numbers from composites~\cite{Delahaye:2013}, and many questions about primes and their distribution amongst the set of natural numbers still remain open. Indeed, most of the knowledge about the sequence of primes stands on unproved theorems and conjectures.

The mystery of primes is not a mere conundrum of pure mathematics. Unexpected connections can be discovered between primes and different topics in physics. For instance, the Riemann zeta function $ \zeta(s) $ --a sum over all integers equivalent to a product over all primes-- has been considered as a partition function~\cite{partition1,partition2,partition3} such that its sequence of non-trivial zeros --encoding information about the sequence of primes-- is similar to the distribution of eigenvalues of random Hermitian matrices used in classically chaotic quantum systems to describe the energy levels in the nuclei of heavy elements~\cite{Forrester:1996}. This idea traces back to the Hilbert-P\'{o}lya conjecture~\cite{hilbert}, which states that the zeros of the $ \zeta(s) $ function might be the eigenvalues of some Hermitian operator on a Hilbert space. Indeed, the Riemann zeta function plays an integral role not only in quantum mechanics but in different branches of physics, from classical mechanics to statistical physics~\cite{Schumayer:2011}. The interpretation of prime numbers or the Riemann zeta zeros as energy eigenvalues of particles appears also in statistical mechanics, as illustrated for instance by the Riemann gas concept as a toy model for certain aspects of string theory~\cite{Spector:1998}. Recently, interesting connections have also been found between primes and self-organized criticality~\cite{Luque:2008fk}, or primes and quantum computation~\cite{shor,latorre} (see~\cite{watkins} for an extensive bibliographical survey between the connection of number theory and physics). The importance of primes transcend theoretical aspects, and practical applications include public key cryptography algorithms~\cite{rsa} and pseudorandom number generators~\cite{Blum:1986fk}. 

One of the most promising approaches to solve the enigmas of number theory is the use of probability theory and stochastic processes. Akin to chaotic dynamical systems, prime numbers, albeit purely deterministic, appear to be scattered throughout natural numbers in a non-homogeneous random fashion. Indeed, for $n\gg 1$ the probability that a randomly chosen number in a ``small'' neighborhood of $n$ is prime is given by~\footnote{\textit{``As a boy I considered the problem of how many primes there are up to a given point. From my computations, I determined that the density of primes around $n$ is about $1/\ln{n}$''}~---~ Carl Friedrich Gauss, 1849}
\begin{equation}
P_n \sim \frac{1}{\ln{n}}.
\label{eq:density}
\end{equation}
This is equivalent to the well-known prime number theorem~\cite{theorem}, which states that the prime counting function $ \pi(N) $ --counting the number of primes up to $N$-- approaches $N/\ln N$ in the limit of $N \to \infty$, i.e.,
\begin{equation}\label{theorem}
\pi(N) \sim \int_2^N \frac{dx}{\ln{x}} \equiv \mbox{Li}(N) \sim \frac{N}{\ln{N}},
\end{equation}
where $\mbox{Li}(N)$ is the offset logarithmic integral function. Taking advantage of this apparent randomness, Cram\'er formulated a simple model~\cite{cramer:1935,cramer:1936} where each integer $n$ is declared as a ``prime'' with independent probability given by Eq.~(\ref{eq:density}). The model --that generates sequences of random primes that are, obviously, in agreement with the prime number theorem-- allowed him to ``prove'', in a probabilistic sense, his famous conjecture about gaps between consecutive primes~\cite{cramer:1936}. 

Cram\'er's probabilistic model plays, still today, a fundamental role when formulating conjectures concerning primes. However, it presents three major drawbacks. 1) It does not ``explain'' the prime number theorem; instead, it is an input of the model. 2) Random primes in the model are totally uncorrelated whereas there are both short and long range correlations in the sequence of real primes. 3) Finally, it says nothing about the relation between prime and composite numbers. In this paper, we combine a complex network approach with the theory of stochastic processes to introduce a parameter-free non-Markovian dynamical model that naturally generates random primes as well as the relation between primes and composite numbers with remarkable accuracy. Our model is in agreement with Eqs.~(\ref{eq:density}) and (\ref{theorem}) and satisfies a modified version of Cram\'er's conjecture about the statistics of gaps between consecutive primes that seems closer to reality than the original Cram\'er's version. Regarding composites, the model helps us to derive the prime factors counting function, giving the probability of distinct prime factors for any integer.

\section{Bipartite network of natural numbers} 

\begin{figure}[t]
\centering
\includegraphics[width=\columnwidth]{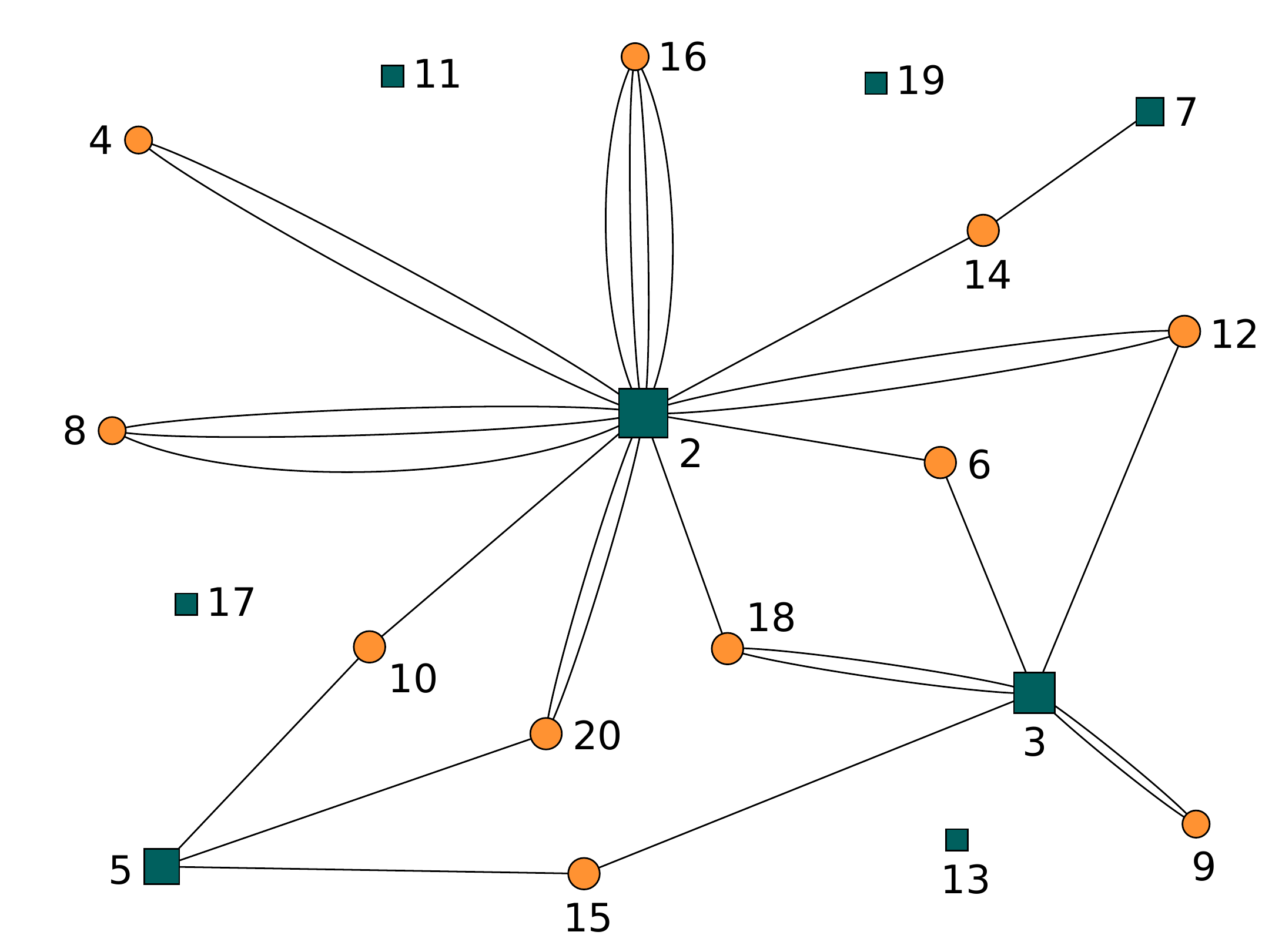}
\caption{\label{primesnet} Example of the bipartite network of natural numbers grown up to size 20. Orange circles represent composite numbers and green squares prime numbers. The degree of a prime, $k_p$, is the number of distinct composites to which it is connected to whereas its strength, $s_p$, is the sum of its weighted connections. Similarly, the degree of a composite, $k_c$, is its number of distinct prime factors and its strength, $s_c$, the total number of prime factors.}
\end{figure}

Primes are the building blocks of natural numbers. The fundamental theorem of arithmetic states that any natural number $ n > 1 $ can be factorized uniquely as
\begin{equation}
\label{eq:arithmetics}
n=p_1^{\alpha_1}p_2^{\alpha_2} \cdots p_k^{\alpha_k}  \cdots
\end{equation}
where $ p_i $ is the $ i $-th prime and $ \alpha_i $ are non-negative integers. From a complex network perspective, natural numbers can be thought of as a weighted bipartite network with two types of nodes, primes and composites. A composite $n$ is linked to primes $p_i$ with weights $\alpha_i$ according to the factorization in Eq.~(\ref{eq:arithmetics}), as shown in Fig.~\ref{primesnet}. 

For a given network size $N$, the probability that a randomly chosen prime inside the network is connected to $k_p$ different composites, that is, the degree distribution $P(k_p)$ for prime numbers, can be exactly determined in terms of the prime counting function as (see Appendix~\ref{appendixA} for details)
\begin{equation}\label{pkbipartite}
P(k_p)=\frac{\pi\left(\frac{N}{k_p+1}\right)-\pi\left(\frac{N}{k_p+2}\right)}{\pi(N)},
\end{equation}
with $k_p=0,1,\cdots, \left\lfloor \frac{N}{2} \right\rfloor$, where $\lfloor x \rfloor$ stands for the floor function. Using the prime number theorem Eq.~(\ref{theorem}), it is easy to see that in the limit $N/k_p \gg 1$ this distribution behaves as $P(k_p)\sim k_p^{-2}$. Quite surprisingly, we obtain a scale-free network with an exponent $-2$, very similar to many real complex networks, like the Internet~\cite{Faloutsos:1999kk}, and similar to the degree distribution of the causal graph of the de Sitter space-time~\cite{Krioukov:2012fk}. As we shall show, this is a consequence of an effective preferential attachment rule induced by the growth mechanism. 

The result in Eq.~(\ref{pkbipartite}) allows us to derive a simple but yet interesting identity relating $\pi(n)$ and the number of distinct prime factors of any integer $n$, $\omega(n)$. We name $\omega(n)$ the prime factors counting function. We start from the trivial identity $[N-1-\pi(N)] \langle k_c \rangle=\pi(N) \langle k_p \rangle$, where $k_c$ is the degree of a composite (or its number of distinct prime factors). Plugging Eq.~(\ref{pkbipartite}) into this identity, we obtain
\begin{equation}
\label{omega0}
\sum_{n=2}^{N} \omega(n)=\sum_{i=1}^{\lfloor N/2 \rfloor} \pi \left( \frac{N}{i}\right).
\end{equation}
Replacing the sum by an integral, we can approximate this expression as
\begin{equation}
\sum_{n=2}^{N} \omega(n)\approx N\int_2^{N} \frac{\pi(x)dx}{x^2}\sim N\ln{\ln{N}}+ \mathcal{O}(N).
\end{equation}
The final asymptotic behavior is directly related to the Hardy-Ramanujan theorem~\cite{hardy-ramanujan}, which now becomes a simple consequence of the prime number theorem. Function $\omega(n)$ can be easily computed from Eq.~(\ref{omega0}) as
\begin{equation}
\omega(n)=\sum_{i=1}^{\lfloor n/2 \rfloor} \left[ \pi \left( \frac{n}{i}\right) -\pi \left( \frac{n-1}{i}\right)\right].
\label{omega}
\end{equation} 
Notice that if $n$ is a composite number, then $\omega(n)$ is, in our network representation, its degree. Therefore, the degree distribution of composite numbers is given by $P(k_c)=\left(\sum_{n=2}^N \delta_{\omega(n),k_c}-\delta_{k_c,1} \pi(N)\right)/(N-1-\pi(N))$. Besides, Eq.~(\ref{omega}) naturally leads to a set of arithmetic functions giving the sum of the prime factors of $n$ raised to any exponent (see Appendix~\ref{appendixC}).

Equations~(\ref{pkbipartite}) and (\ref{omega}) are a remarkable result. Beyond potential applications to find better estimates of function $\omega(n)$, they state that the local properties of the network of natural numbers are fully determined by the prime counting function $\pi(N)$ alone. We then expect that any model producing random versions of the network that is able to reproduce well the prime counting function, $\pi(N)$, will also reproduce well the large scale of the real network topology. 

\begin{figure}[t]
\centering
\includegraphics[width=\linewidth]{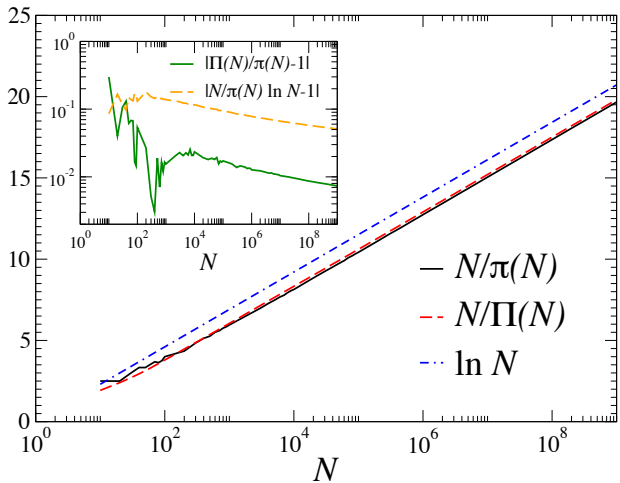}
\caption{Comparison of the prime counting function $\pi(N)$, the prime number theorem Eq.~(\ref{theorem}), and the prime counting function of our random model $\Pi(N)$, averaged over $1000$ realizations. The inset shows the corresponding relative errors. The relative error of the random model is one order of magnitude smaller than the one of Eq.~(\ref{theorem}).\label{fig2}}
\end{figure}

\section{Modeling the evolution and structure of natural numbers} 

The order relation implicit in the natural numbers allows us to consider the bipartite network representation of natural numbers as a growing system. In the growing process, natural numbers join the network sequentially and try to connect to already existing primes. Those new numbers that succeed in this process are said to be composites, otherwise, they become prime numbers. In this paper, we show that a very simple connection rule based upon a soft version of Eq.~(\ref{eq:arithmetics}) generates networks with the same architecture as that of the real network of natural numbers. Taking advantage of the apparent randomness of prime numbers, we develop a stochastic model that generates growing bipartite natural number networks connecting random primes with composites. The growth process only imposes two basic facts trivially implied by the fundamental theorem of arithmetic, that is, that the product of the prime factors of a natural number $n$ must be $n$, and that $n$ can have no more than one prime factor larger than $\sqrt{n}$. The model starts by assuming that number $2$ is a prime and adds natural numbers $n \ge 3$ sequentially. It proceeds as follows 
\begin{enumerate}
\item
Each new number $n$ that joins the network tries to connect to already existing random primes $p_i \le \sqrt{n}$ with independent probabilities $1/p_i$ one by one, starting from the smallest prime, until the first connection is stablished.
\item
If number $n$ first connects to an existing prime $p$ smaller or equal to $\sqrt{n}$, it keeps trying to connect sequentially to existing primes in the range $[R_{m}, R_{M}]$, with $R_{m}=p$ and $R_{M}=\sqrt{n'}$, and $n'=\frac{n}{p}$. Each time $n$ connects to a new random prime $p'$ the range is redefined with $R_{m,new}=p'$ and $n'_{new}=\frac{n'_{old}}{p'}$. If $p' > R_{M,new}$ or $n$ does not get new connections in the evaluation range, $n$ is connected to the prime closest to $R_{M}^2$ and a new node $n+1$ is added to the system.
\item
If number $n$ does not connect to any existing prime smaller or equal to $\sqrt{n}$, it is declared as a prime and a new number $n+1$ is added to the system.
\end{enumerate}
The intuition behind the second step in our model is as follows. In the case of the real primes, a composite number $n$ must have at least a prime factor smaller or equal to $\sqrt{n}$. Let $p$ be the smallest prime factor of $n$. Then, $n/p$ is also an integer number that is either a prime or, else, it can be expressed as a product of prime factors. However, in the latter case the smallest prime factor of $n/p$ cannot be smaller than $p$ because this  would contradict the assumption that $p$ is the smallest prime factor of $n$. Then, the smallest prime factor of $n/p$, let it be $p'$, must lie in the closed interval $[p,\sqrt{n/p}]$. The same logic can now be applied to the prime factors of the ratio $n/(pp')$ until $n$ is fully factorized. Our model tries to mimic in a stochastic manner this factorization property of composite numbers, with the difference that, in our case, $n/p$ may not be an integer. Thus, at the end of a stochastic realization of our model, every number $n$ is either declared as a prime or it is a composite such that the product of its prime factors is approximately $n$.
\begin{figure}[t]
\centering
\includegraphics[width=\linewidth]{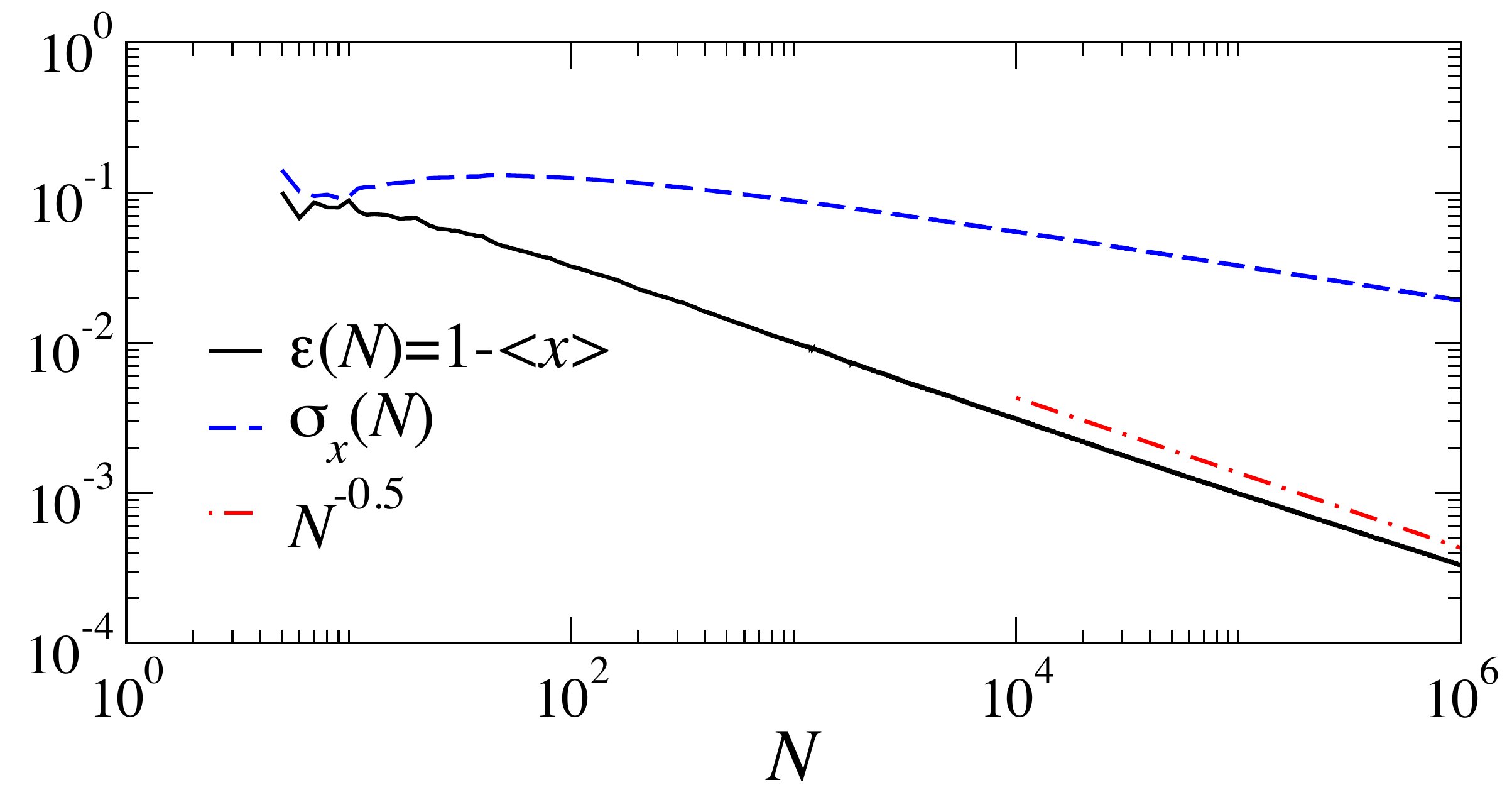}
\caption{Average relative error $\epsilon(N)=1-\langle x \rangle$ between a composite number and its factorization in the network as a function of the system size $N$ and the standard deviation of the ratio $x$, $\sigma_x(N)$. \label{fig3}}
\end{figure}

It is worth noticing the following properties of the model. i) The model has no tunable parameters. ii) It is a generative model, in the sense that the model generates simultaneously the number of primes and how primes and composites are connected. iii) The model is able to generate multiple connections between composite and a prime numbers with no extra mechanism. iv) The model is non-Markovian because the probability of a number being a prime depends on the whole history of the stochastic process. At this respect, it is important to notice that all results in this paper are considered to be averages over all histories of the stochastic process. We also notice that the first step of the algorithm is similar to the random sieve proposed by Hawkins~\cite{hawkins:1957,hawkins:1974,Bui:2006,Lorch:2007}. The main difference being that the random sieve does not provide connections between composite and prime numbers.

\subsection{The prime counting function} 

The analytical treatment of the model is quite involved due to its non-Markovian character (see Appendix~\ref{appendixD}). However, it is possible to work out a relatively simple mean field approximation. For instance, the probability that number $n$ is a prime according to the model, $P_n$, satisfies the following recurrence relation
\begin{equation}
P_n=e^{\displaystyle{\sum_{i=2}^{\sqrt{n}}\ln{\left[1-\frac{P_i}{i}\right]}}} \approx e^{-\displaystyle{\int^{\sqrt{n}} \frac{P_x}{x} dx}},\\
\label{recurrence1}
\end{equation}
where in the last term we have considered $n$ as a continuous variable and approximated $\ln{\left[1-\frac{P_i}{i}\right]}$ by $-\frac{P_i}{i}$. It is easy to see that Eq.~(\ref{recurrence1}) is equivalent to the following non-linear and non-local differential equation
\begin{equation}
\frac{dP_n}{d n}=-\frac{P_n P_{\sqrt{n}}}{2n}.
\end{equation}
Although the full analytical solution of this equation is difficult to obtain, it is quite easy to check that, asymptotically, $P_n$ behaves as $P_n \sim 1/\ln{n}$ and, thus, our model satisfies the prime number theorem as an emerging property. Figure~\ref{fig2} shows a comparison between the real $\pi(N)$, the one generated by our model $\Pi(N)$, and Eq.~(\ref{theorem}). As expected, $\lim_{N\rightarrow \infty} \pi(N)/\Pi(N)=1$. However, for finite sizes the relative error of our model with respect to the real $\pi(N)$ is one order of magnitude smaller than the one given by Eq.~(\ref{theorem}).
\begin{figure}[t]
\centering
\includegraphics[width=\linewidth]{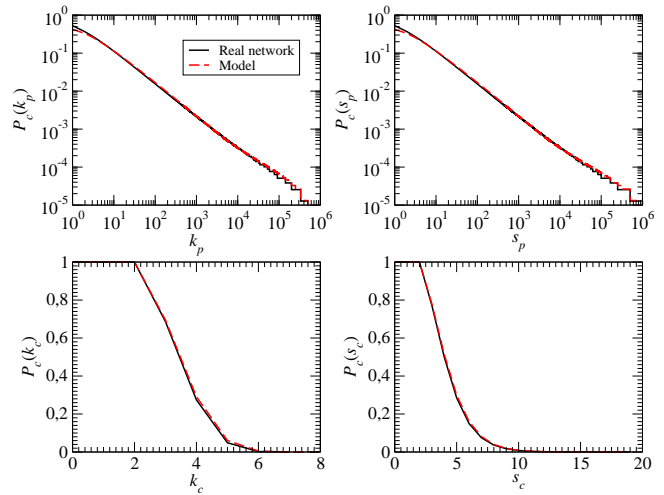}
\caption{Comparison between the complementary cumulative distribution functions of the real bipartite network of natural numbers of size $N=10^6$ and the network generated by our model averaged over $1000$ realizations. The left column shows the unweighted properties and the right column the weighted ones. The legend explaining line types applies to the four plots. \label{fig4}}
\end{figure}

\subsection{Network properties} 

One of the strengths of our model lays in its ability to reproduce, not only the sequence of primes, but also the connections of each composite number. To check to what extent our model fulfills the fundamental theorem of arithmetic, we measure the relative error between a composite and its factorization according to the model, $\epsilon(N)$, as follows. Let $c_i$ be the $i$th composite in a network of size $N$ and let $\bar{c}_i$ be its factorization, then we define $x_i \equiv \bar{c}_i/c_i $. The relative error is then $\epsilon(N) \equiv 1-\langle x \rangle=1-(N-1-\Pi(N))^{-1} \sum_{i}\bar{c}_i/c_i$, where $\langle \cdot \rangle$ means the population average. In Fig.~\ref{fig3}, we show $\epsilon(N)$ as a function of the system size $N$ averaged over $1000$ network realizations. As it can be seen, this error decreases as a power law of the size of the system $\epsilon(N) \sim N^{-\alpha}$ with $\alpha \approx 0.5$. We also show the standard deviation of $x_i$, which also approaches zero in the large system size limit. These two results indicate that the model fulfills the fundamental theorem of arithmetic for relatively small numbers with high accuracy. 

The model also does an excellent job at reproducing well the large-scale topology of the real network. The left column in Fig.~\ref{fig4} shows the complementary cumulative degree distributions of primes and composites as compared to the real ones for the network grown up to $N=10^6$. In both cases the agreement is excellent. The right column in Fig.~\ref{fig4} shows the strengths distributions for primes and composites, that is, the equivalent to the left column measures when multiple links between primes and composites are considered (see Fig.~\ref{primesnet}). Again, the agreement between the model and the real network is excellent. This result is particularly interesting as it shows that our model is able to capture statistical properties of the multiplicities of  composites' factorizations, i. e. the $\alpha$s in Eq.~(\ref{eq:arithmetics}). In particular, it recovers that $P_{c}(s_{p})$ behaves asymptotically as $s_{p}^{-2}$, as expected from the almost linear correlation between strength and degree. Other topological properties are explored in Appendices~\ref{appendixB} and \ref{appendixE}. For instance, it is possible to show that the model satisfies the Erd\"os-Kac theorem~\cite{erdos-kac}, which states that $(\omega(n)-\ln{\ln{n}})/\sqrt{\ln{\ln{n}}}$ is, {\it de facto}, a random variable that follows the standard normal distribution.

\subsection{The Cram\'er's conjecture revisited}

Cram\'er's conjecture provides an absolute upper bound on the gaps between consecutive primes. Using his model, Cram\'er was able to prove that~\cite{cramer:1936} 
\begin{equation}
\limsup_{i\rightarrow \infty }\frac{p_{i+1}-p_i}{\ln^2{p_i}}=1
\end{equation}
and conjectured that the same relation also holds for real primes. Here, we study the statistics of prime gaps in our model and refine Cram\'er's conjecture for real primes. We start by noticing that in our model, all numbers between two perfect squares have the same probability of being primes and, more importantly, they are conditionally independent given their common history. Therefore, as a first approximation, we consider that every number in the interval $[m^2,(m+1)^2)$; $m=2,3,\cdots$ has an independent probability $P_n=1/\ln{n}$ of being a prime, where $n=m^2$. Under this assumption, the probability that a given gap $G$ within the interval is smaller than $g$ is $\mbox{Prob}\{ G< g\}=1-(1-P_n)^{g-1}$~\footnote{We are implicitly assuming that the length of the interval is much larger than the typical gap and, therefore, that the number of gaps within the interval is large.}. If we assume that there are $N_{G}=2 \sqrt{n} P_n$ gaps within the interval, the probability that the largest gap $G_{m}$ within the interval is smaller than $g_m$ is 
\begin{equation}
\mbox{Prob}\{ G_{m}<g_m\}=[1-(1-P_n)^{g_m-1}]^{N_{G}}.
\end{equation}
The average largest gap can be evaluated from this expression, yielding
\begin{equation}
\langle G_{m} \rangle=\left(\frac{1}{P_n}-\frac{1}{2}\right) H_{N_{G}}+\mathcal{O}(P_n) \sim \frac{1}{2}\ln^2 n,
\label{eq:averagegapmax}
\end{equation}
where $H_{N_{G}}=\sum_{k=1}^{N_{G}}k^{-1}$ is the harmonic number (interestingly, a similar approach has been recently proposed in~\cite{Wolf:2014fk}). We can now define the normalized largest gap as $\overline{G}_m \equiv G_{m}/\langle G_{m} \rangle$, which distribution function satisfies
\begin{equation}
\mbox{Prob}\{ \overline{G}_m<\overline{g}_m\} \sim e^{-N_{G}^{1-\overline{g}_m}}.
\label{eq:gapdistribution}
\end{equation}
In the limit $n \rightarrow \infty$, $N_G  \rightarrow \infty$ and this distribution becomes a step function (although very slowly). Thus, the largest gap stops being a random variable to become a deterministic quantity equal to $\ln^2 n/2$. Notice that this bound is twice as small as the bound given by Cram\'er's conjecture, apparently suggesting that it could be false for real primes. 
\begin{figure}[t]
\centering
\includegraphics[width=\linewidth]{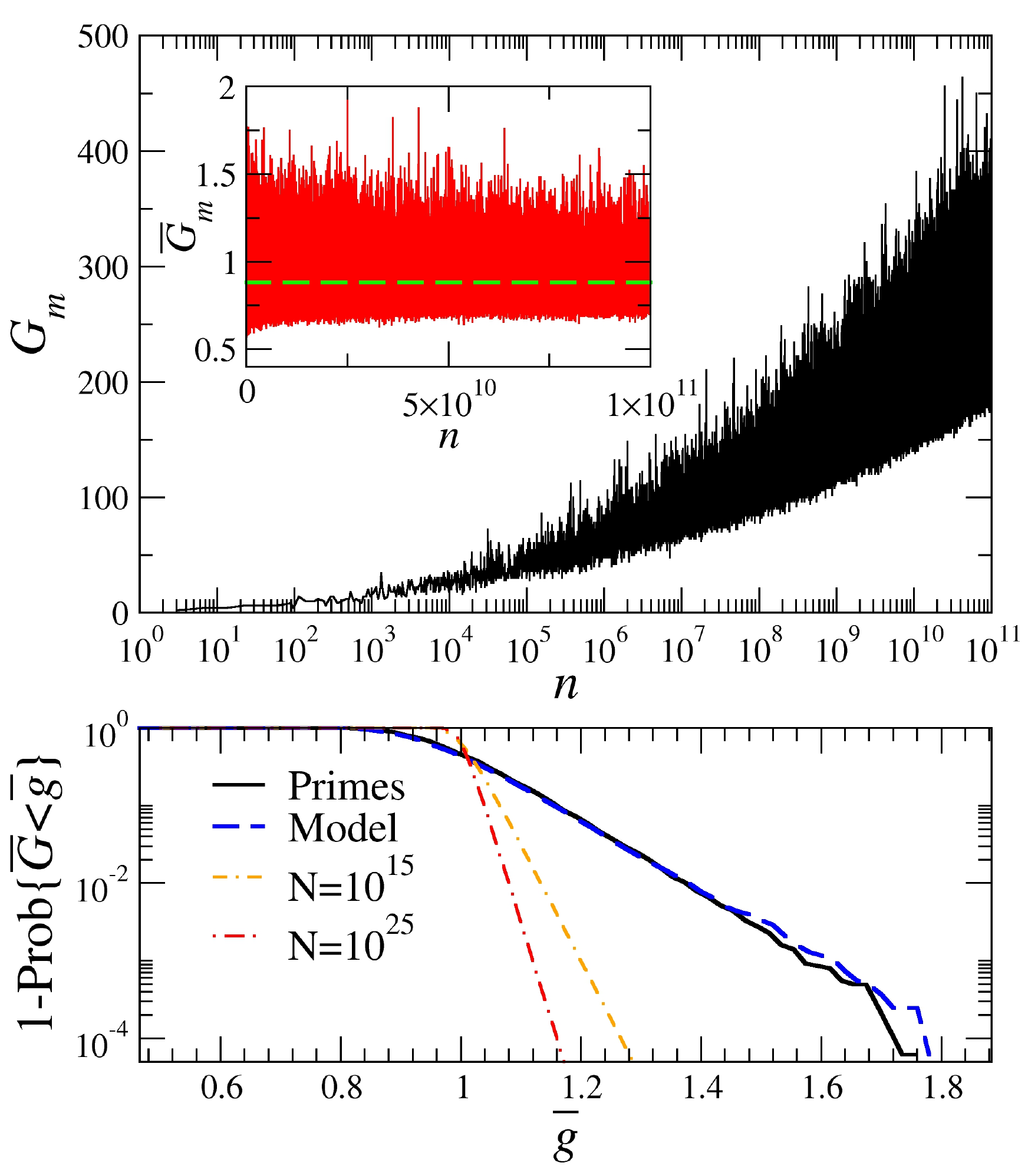}
\caption{Gaps between primes. Top. Series of largest gaps between real primes in intervals between perfect squares. Top Inset. The same series normalized by using Eq.~(\ref{eq:averagegapmax}). In both plots, primes are considered up to $10^{11}$. Bottom. Complementary cumulative distribution function of the normalized largest gaps for real primes and the model in the range $[9\times 10^{10},10^{11}]$. To make evident the slow convergence of the distribution, we also show extrapolations from Eq.~(\ref{eq:gapdistribution}) for $N=10^{15}$ and $N=10^{25}$. \label{fig5}}
\end{figure}
To check our prediction, we compute the gaps between real primes up to $10^{11}$. We divide this set in intervals between perfect squares and for each such interval we evaluate the largest gap. The top plot in Fig.~\ref{fig5} shows the series of largest gaps and the inset shows the normalized largest gaps obtained by using Eq.~(\ref{eq:averagegapmax}). As it can be seen, after normalization, the series becomes a stationary one but its average is not 1, as we would expect from our model, but $2 c \approx 0.88$, with $c$ a constant below $1/2$. As we see, our model suffers from the same problems affecting Cram\'er's model in what respect short range correlations induced by small primes. For instance, the probability of $n$ being a prime if $n-1$ is a prime is zero for real primes whereas our model would predict a non-zero probability; in addition, the probabilistic prediction that the number of primes in a short interval of length $y$ about $x$ is given by $y/\ln x $ was proved false by Maier~\citep{Maier:1985,Granville:1994}. Some other deviations from real primes on a very large scale have also been reported~\citep{Odlyzko:1985,Rubinstein:1994}. In the case of Cram\'er's model, it is possible to make heuristic corrections allowing one to reach right answers on several properties of real primes, like the number of twin primes below $N$~\cite{Hardy:1923}. In general, these corrections have only a numerical effect on the studied property since the bear model already predicts the right asymptotic behavior as a function of $N$. The same type of heuristics can be, in principle, applied to our model and we expect them to account for the observed discrepancy. For instance, a simple modification assumes that the probability of $n$ being a prime is zero if the previous number is a prime whereas it is $(\ln{n}-1)^{-1}$ otherwise. This simple modification preserves the prime number theorem and leads to a better estimate of constant $2c \approx 0.92$.

\begin{figure}[t]
\centering
\includegraphics[width=\linewidth]{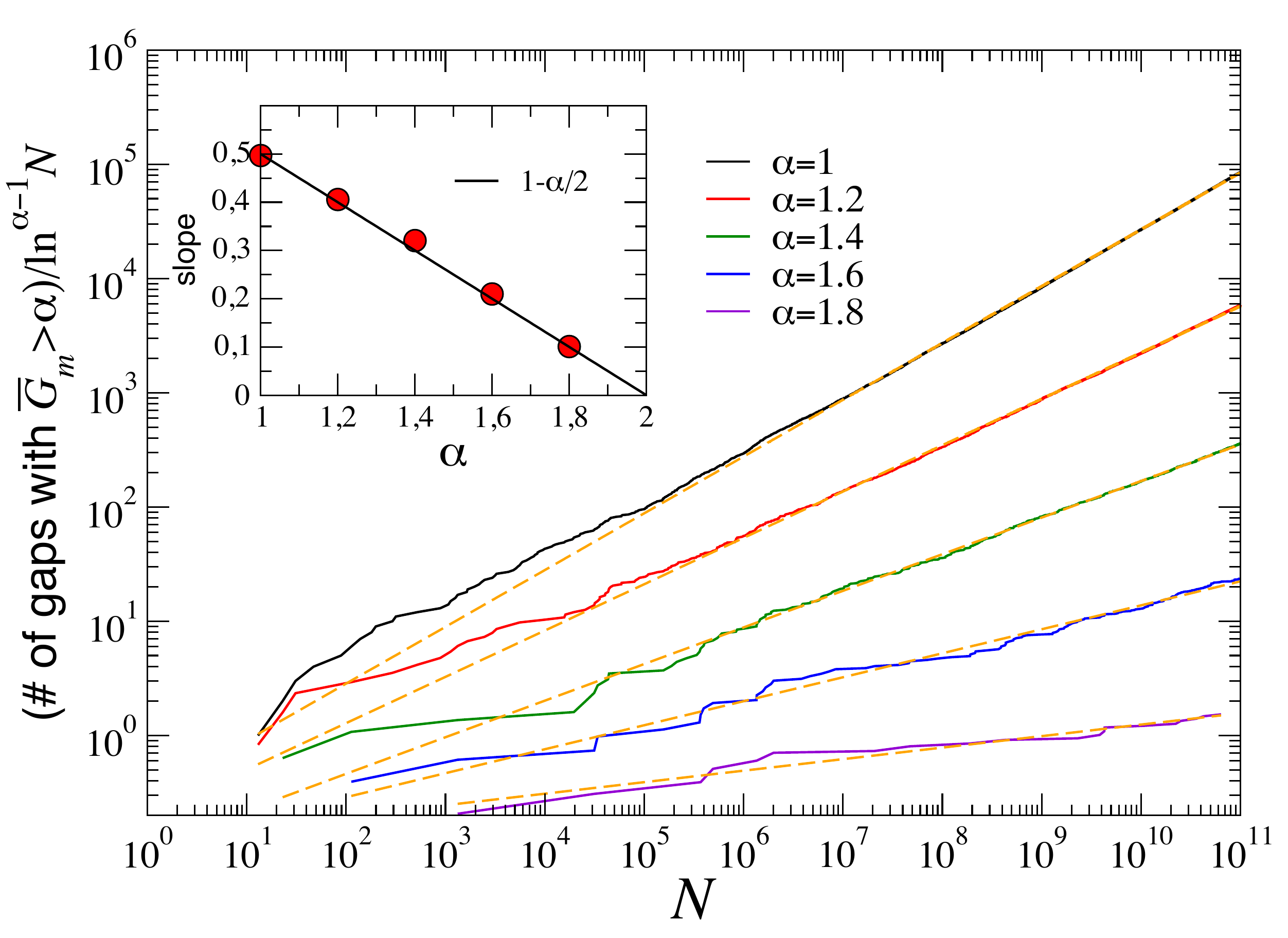}
\caption{Number of gaps with $\overline{G}_m>\alpha$ for different values of $\alpha$ as a function of $N$ re-scaled by the factor $\ln^{\alpha-1}N$. According our estimates, this should behave as a power law of the form $N^{1-\alpha/2}$. Dashed lines are power law fits which exponents are shown in the inset plot and compared to the theoretical prediction $1-\alpha/2$. \label{fig6}}
\end{figure}

Even more interesting is the analysis of the fluctuations of the normalized largest gaps around their average. 
A preliminary analysis of their distribution suggests that largest gaps of real primes behave as in the model after a global rescaling. Thus, to have a coherent comparison between the model and real primes, we divide the series shown in the inset of Fig.~\ref{fig5} by $2 c$ so that its average is equal to 1, like in the model. We then evaluate the complementary cumulative distribution function for all largest gaps in the range $[9\times 10^{10},10^{11}]$ and compare it with the one obtained from numerical simulations of our model, see bottom plot in Fig.~\ref{fig5}. Interestingly, both distributions are nearly indistinguishable. This implies that fluctuations of largest gaps for real primes are governed asymptotically by the distribution Eq.~(\ref{eq:gapdistribution}). From this equation, we can evaluate the expected number of gaps up to $N$ that are above a certain fraction $\alpha$ of the average largest gap, with $\alpha \ge 1$, that is,
\begin{equation}
\mbox{\# gaps with } \overline{G}_m>\alpha  \approx \sum_{n=1}^{\sqrt{N}} \left( \frac{\ln{n}}{n} \right)^{\alpha-1}.
\end{equation}
This quantity diverges when $1 \le \alpha < 2$ as $\mathcal{O}(N^{1-\alpha/2} \ln^{\alpha-1}N)$ and as $\mathcal{O}(\ln^2 N)$ for $\alpha=2$. Putting all the pieces together, we refine Cram\'er's conjecture as follows. For all real prime gaps $G_i \equiv p_{i+1}-p_i$, with $p_i<N$ and $N\rightarrow \infty$, we have
\begin{equation}
\hspace{-0.14cm}
\begin{array}{ll}
G_i <  \alpha c \ln^2{p_i} & \mbox{ for all but $\mathcal{O}(N^{1-\frac{\alpha}{2}} \ln^{\alpha-1}N)$ gaps}\\
G_i<  2c \ln^2{p_i} &\mbox{ for all but $\mathcal{O}(\ln^2 N)$ gaps} 
\end{array}
\label{cramer_improved}
\end{equation}
For any $\alpha>2$, the number of gaps above this threshold is $\mathcal{O}(1)$. Notice however that this asymptotic behavior is only reached for extremely large values of $N$. For not so large values it is better to replace $\ln^2 p_i$ in Eq.~(\ref{cramer_improved}) by $2 \left[\ln{p_i}-1/2\right]\left[\ln{\left( 2 \sqrt{p_i}/\ln{p_i}\right)}+\gamma \right]$, with $\gamma$ the Euler-Mascheroni constant, as derived from Eq.~(\ref{eq:averagegapmax}). We check these predictions for all gaps up to $10^{11}$ in Fig.~\ref{fig6}. We measure empirically the number of gaps that, up to a given size $N$, satisfy $G_i>2 \alpha c \left[\ln{p_i}-1/2\right]\left[\ln{\left( 2 \sqrt{p_i}/\ln{p_i}\right)}+\gamma \right]$ and compare them with the predictions in Eq.~(\ref{cramer_improved}). Aside from statistical errors, our predictions agree well with the empirical measures. 

\section{Conclusions}

Probabilistic approaches to understand usual patterns of primes as well as their extreme statistics brought a new perspective to the study of prime numbers. The big first step by Cram\'er was significantly developed afterwards bringing this kind of approach to maturity. With our work, we introduce a new dimension that allows us to understand primes and their statistical properties not in isolation but as building blocks of natural numbers.  We have introduced a parameter-free non-Markovian stochastic model based on a bipartite complex network representation that naturally generates random primes as well as the relation between primes and composite numbers with remarkable accuracy. Our model satisfies the Erd\"{o}s-Kac theorem, as well as the prime number theorem and a refined version of Cram\'er's conjecture about the statistics of gaps between consecutive primes that seems closer to reality than the original Cram\'er's version. Even though we are still unable to fully understand the finer details about primes and the complex architecture of natural numbers, probabilistic models like ours provide valuable tools helping to elaborate conjectures about primes and, perhaps, also to prove results on number theory. Beyond the implications in mathematics, our stochastic model generates the sequence of random primes and some of their statistical correlations as an emergent property, which allows probabilistic computations of number theoretical approaches to open problems in physics involving the Riemann zeta function, which plays an integral role in different branches from quantum mechanics to condensed matter.

\begin{acknowledgments}
We thank Dmitri Krioukov for useful comments and suggestions. We acknowledge support from the James S. McDonnell Foundation 21st Century Science Initiative in Studying Complex Systems Ð Scholar Award; the ICREA Academia prize, funded by the {\it Generalitat de Catalunya}; MICINN projects No.\ FIS2010-21781-C02-02 and BFU2010-21847-C02-02; {\it Generalitat de Catalunya} grant No. 2014SGR608; and the Ram\'on y Cajal program of the Spanish Ministry of Science.
\end{acknowledgments}

\onecolumngrid
\appendix

\section{Bipartite network representation of natural numbers}
\label{appendixA}
In this section we derive the expressions that characterise the bipartite network representation of natural numbers presented in the paper.

\subsection{Degree distribution}
The degree distribution for primes in the network can be derived reasoning as follows: a prime number $p > N/2$ has degree $k_{p}(p)=0$ since its product by any other prime number is greater than $N$ and, hence, it cannot belong to the network (the subscript in $k_{p}$ is used to denote the degree of primes; we use $k_{c}$ to refer to the degree of composites). Identically, if $N/3 < p \leqslant N/2$, $p$ has a multiple which belongs to the network $(2p \leqslant N)$. In general,
\begin{equation}\label{interval}
p \in \left( \frac{N}{n+1}, \frac{N}{n} \right] \Leftrightarrow k_{p}(p) = n - 1,
\end{equation}
since $mp \leqslant N, ~ m=2,\ldots,n$ but $(n+1)p > N$. This directly leads to the expression for $P(k_{p})$,
\begin{equation}\label{pk_p}
P(k_{p}) = \frac{\# \{ p : p ~ \text{prime} : k_{p}(p)=k_{p} \}}{\# \{ p : p ~ \text{prime} \leqslant N \}} = \frac{\# \{ p : p ~ \text{prime} \in \left( \frac{N}{k_{p}+2}, \frac{N}{k_{p}+1} \right] \}}{\# \{ p : p ~ \text{prime} \leqslant N \}} = \frac{\pi \left( \frac{N}{k_{p}+1} \right) - \pi \left( \frac{N}{k_{p}+2} \right)}{\pi \left( N \right)}.
\end{equation}
This expression is, interestingly, similar to a probability measure with multifractal properties used in~\cite{wolf:1989}.
We can derive an approximation for Eq.~\eqref{pk_p} using the fact that, according to the prime number theorem,
\begin{equation}\label{pnt}
\lim_{x \to \infty} \frac{\pi(x)}{x/\ln(x)} = 1,
\end{equation}
We first evaluate the complementary cumulative distribution function $P_c(k_p)=\sum_{k=k_p} P(k_p)$, which reads
\begin{equation}
P_c(k_p)=\frac{\pi \left( \frac{N}{k_{p}+1} \right)}{\pi(N)}.
\end{equation}
Using the prime number theorem, in the limit $N/k_{p} \gg 1$ this function behaves as
\begin{equation}
P_c(k_p) \approx \frac{1}{k_p (1-\frac{\ln{k_p}}{\ln{N}})} \sim \frac{1}{k_p},
\end{equation}
from where it follows that the degree distribution behaves nearly as a power law
\begin{equation}\label{pk_p_powerlaw}
P(k_{p}) \sim k_{p}^{-2}.
\end{equation}

Another useful relation is
\begin{equation}\label{degree}
k_{p}(p) = \left\lfloor \frac{N}{p} \right\rfloor - 1,
\end{equation}
which can be proved considering \eqref{interval}
\begin{equation*}
p \in \left( \frac{N}{n+1}, \frac{N}{n} \right] \Leftrightarrow \frac{N}{p} \in \left[ n, n+1 \right) \Leftrightarrow \left\lfloor \frac{N}{p} \right\rfloor = n \Leftrightarrow k_{p}(p) = \left\lfloor \frac{N}{p} \right\rfloor - 1.
\end{equation*}

\subsection{Strength of a prime number}
The expression for the strength of a prime number $p$ in the network of size $N$ is
\begin{equation}\label{strength}
s_{p}(p) = \sum \limits_{n=1}^{\left\lfloor \log_{p} N \right\rfloor} \left\lfloor \frac{N}{p^{n}} \right\rfloor - 1.
\end{equation}
The explanation of this formula is rather straightforward. The prime $p$ inside the bipartite network is connected to $\left\lfloor N/p \right\rfloor - 1$ composites (Eq.\eqref{degree}). Nevertheless, $\left\lfloor N/p^2 \right\rfloor$ of these composites can be divided by $p$ twice. In general, there are $\left\lfloor N/p^n \right\rfloor$ composites which can be divided by $p$ $n$ times. Since the strength of the prime $p$ is defined as the sum of the weights of all its connections, we can simply sum all these terms as
\begin{equation*}
s_{p}(p) = k_{p}(p) + \sum \limits_{n=2}^{\infty} \left\lfloor \frac{N}{p^{n}} \right\rfloor = \sum \limits_{n=1}^{\infty} \left\lfloor \frac{N}{p^{n}} \right\rfloor - 1.
\end{equation*}
An upper limit for the sum can be found by taking into account the fact that, if $p^n > N \Rightarrow N/p^n < 1$ and, hence, such term does not contibute to the sum. Let us then find the values of $n$ which need to be considered,
\begin{equation*}
\left\lfloor \frac{N}{p^n} \right\rfloor > 0 \Leftrightarrow \frac{N}{p^n} \geqslant 1 \Leftrightarrow p^n \leqslant N \Leftrightarrow n \leqslant \log_{p} N.
\end{equation*}
This allows us to write the upper limit in Eq.\eqref{strength}, since the last term to be added is the one for $n=\left\lfloor \log_{p} N \right\rfloor$.

\subsection{Strength distribution}
A reasonable approximation of the strength as a function of the degree $k_{p}$ is given by
\begin{equation}\label{sk}
s_{p}(k_{p}) \sim \frac{N \left(k_{p}+1 \right)}{N - \left(k_{p}+1 \right)} - 1,
\end{equation}
which shows that weights do not play an important role in our representation since, for small values of $k_{p}$, Eq.~\eqref{sk} exhibits a linear behaviour ($s_{p}(k_{p}) \sim k_{p}$). This result is a consequence of the fact that only primes less or equal to $\sqrt{N}$ have connections with weight greater than 1, which implies that the fraction of nodes for which this is possible, $1/\sqrt{N}$, tends to zero in the thermodynamic limit. Eq.~\eqref{sk} can be derived by approximating Eq.\eqref{strength} as
\begin{equation}\label{step1}
s_{p}(p) = \sum \limits_{n=1}^{\left\lfloor \log_{p} N \right\rfloor} \left\lfloor \frac{N}{p^{n}} \right\rfloor - 1 \sim \sum \limits_{n=1}^{\infty} \frac{N}{p^n} - 1 = \frac{N}{p-1} - 1.
\end{equation}
We can finally use Eq.~\eqref{degree} to give an approximate value of $p(k_{p})$, i.e. a prime with degree $k_{p}$, 
\begin{equation}\label{inversion}
k_{p}(p) = \left\lfloor \frac{N}{p} \right\rfloor -1 \Rightarrow p \sim \frac{N}{k_{p}+1}.
\end{equation}
The substitution of Eq.~\eqref{inversion} into Eq.~\eqref{step1} yields Eq.~\eqref{sk}.

The cumulative strength distribution can also be derived as follows. From Eq.~(\ref{step1}), we see that any prime $p$ such that
\begin{equation*}
p \gtrsim \frac{N}{s_{p}+1}+1
\end{equation*}
must have strength less or equal to $s_{p}$. We can therefore approximate $P_{c}(s_{p}) = \text{Prob} \left\{ S > s_{p} \right\} = 1 - \text{Prob} \left\{ S \leq s_{p} \right\}$, where $S$ stands for the strength of a randomly chosen prime, as
\begin{align*}
P_{c}(s_{p}) &\sim 1 - \frac{\pi \left( N \right) - \pi \left( \frac{N}{s_{p}+1}+1 \right)}{\pi \left( N \right)} = \frac{\pi \left( \frac{N}{s_{p}+1}+1 \right)}{\pi \left( N \right)} \sim \frac{\pi \left( \frac{N}{s_{p}+1} \right)}{\pi \left( N \right)} \\
&\sim \frac{N}{\left( s_{p} + 1 \right) \ln \left(\frac{N}{s_{p}+1}\right)} \frac{\ln N}{N} = \frac{1}{1 - \frac{\ln (s_{p}+1)}{\ln N}} \frac{1}{s_{p} + 1} \sim s_{p}^{-1},
\end{align*}
so we see that, indeed, $P(s_{p}) \sim s_{p}^{-2}$.

\section{One-mode projection}
\label{appendixB}
Given a bipartite network, we can build a new graph composed exclusively of nodes belonging to one of its classes by performing the so called one-mode projection. Since no pair of these nodes can be initially connected by the definition of bipartite network, linking must be ruled by some other criteria in the new graph. The most usual one is to establish a connection between two nodes with a weight equal to the number of common nodes to which they were both connected in the original network. Hence, whenever two nodes had no common neighbours in the bipartite network, they are left unconnected.

\begin{figure}[h!]
\centering
\includegraphics[scale=0.7]{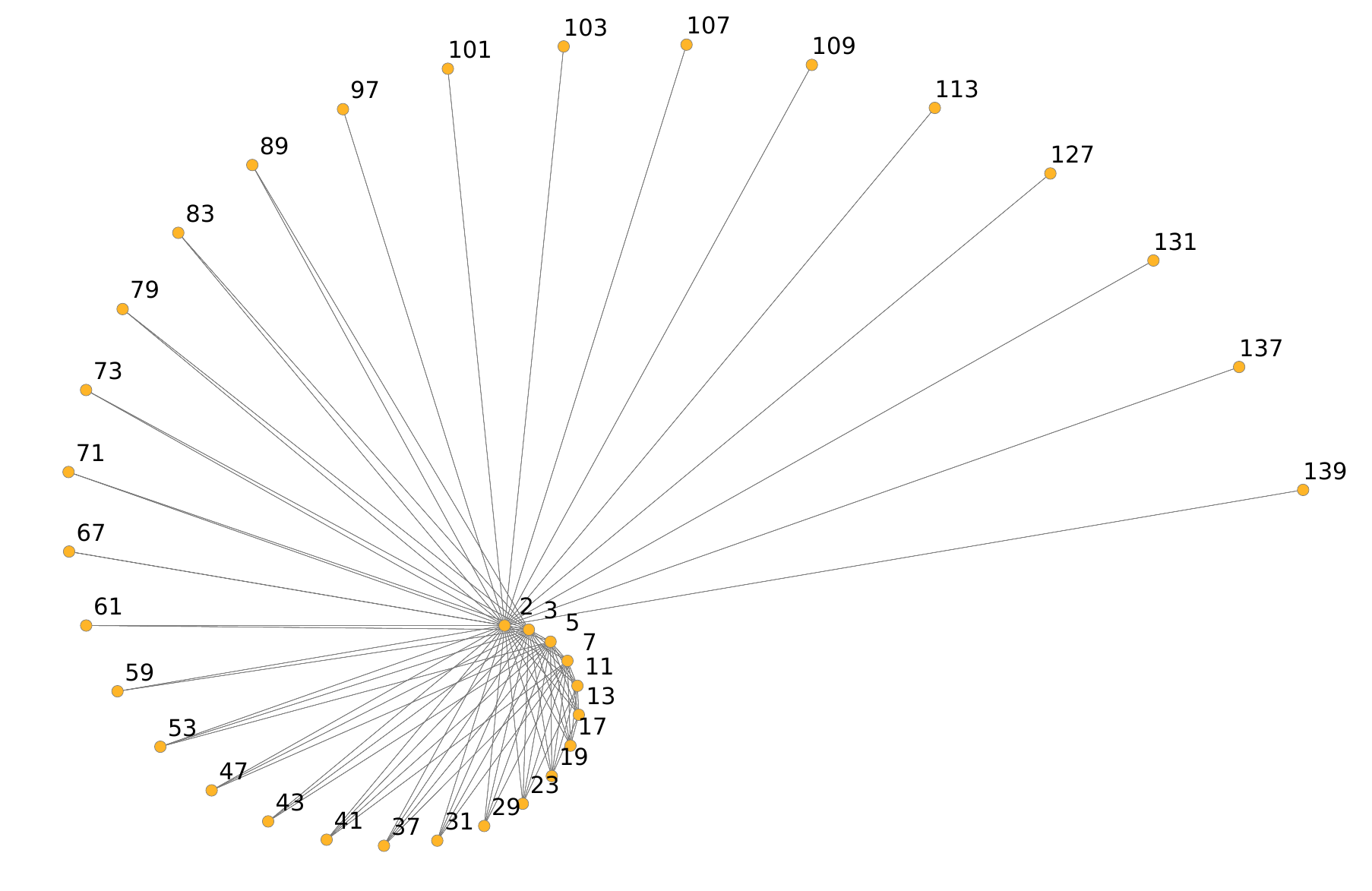}
\caption{\label{omp} One-mode projection for $N=289$. Primes greater than $N/2$ do not appear in the picture since they are all unconnected. Self-loops are not depicted either. Primes $\left\{ 2,3,5,7,11,13,17 \right\}$ form a clique, and there are no links between nodes not belonging to it.}
\end{figure}

In order to deepen into the study of the statistical properties of prime numbers, we have performed a one-mode projection onto that class in the bipartite network discussed so far following the latter criteria (see Fig.~\ref{omp}) and, in addition, allowing self-loops to exist in the resulting graph (whenever a perfect power of a prime exists in the bipartite network, we regard that prime as connected to itself, thus forming a self-loop).

As can be seen in Fig.~\ref{omp}, and as the results presented in this section imply, this graph has a structure made of a maximally connected core containing all the primes less or equal to $\sqrt{N}$ that is surrounded by nodes connected to some but not all of the inner nodes. In addition, the inner two prime numbers are, the strongest the connection amongst them. This suggests that this network could exhibit a self-similar behaviour, i.e. it could be statistically invariant under a network renormalization procedure. This interesting property would allow us to predict some of its statistical properties on any scale.

\subsection{Degree distribution}
The degree $k$ of a prime number $p$ in the one-mode projection of a bipartite network of size $N$ is given by
\begin{equation}\label{om_degree}
k(p) = \pi \left( \frac{N}{p} \right).
\end{equation}
This expression is justified as follows: $p$ can be connected to any prime number $p'$ as long as $pp' \leqslant N$. As a consequence, in order to obtain the number of primes $p'$ to which $p$ can be connected, we must count the number of primes $p' \leq N/p$, which is precisely the result in Eq.~\eqref{om_degree}. Notice that, if $p \leq \sqrt{N}$, $p$ is counted as well; hence, this expression takes self-loops into account.

Using $p_{k}$ to denote the $k$-th prime, the degree distribution $P(k)$ is exactly determined by
\begin{equation}\label{pkonemode}
P(k) = \frac{\pi \left( \frac{N}{p_{k}} \right) - \pi \left( \frac{N}{p_{k+1}} \right)}{\pi \left( N \right)}, \quad p_{0} \equiv 1.
\end{equation}
This result starts with the observation that if a prime $p$ has degree $k$, it must be connected to the first $k$ prime numbers $p_1,p_2,...,p_{k}$. Hence, $p p_{k} \leqslant N$ but $p p_{k+1} > N$. In order to count how many primes are subject to these conditions, we must count the number of primes in the interval $p \in \left( N/p_{k+1}, N/p_{k} \right]$, which can be written in terms of th prime counting function as $\pi \left( N/p_{k} \right) - \pi \left( N/p_{k+1} \right)$. Dividing that quantity by the amount of primes in the graph $\pi(N)$ yields Eq.~\eqref{pkonemode}. We must take into account that, in the particular case of $k=0$, we are considering the primes $p$ for which $p p_1 > N$ and $p \leqslant N$, i.e. the primes $p \in \left( N/p_{1}, N \right]$. Defining $p_0 \equiv 1$, the latter equation is extended to that case.

\begin{figure}[h!]
\centering
\includegraphics[scale=0.4]{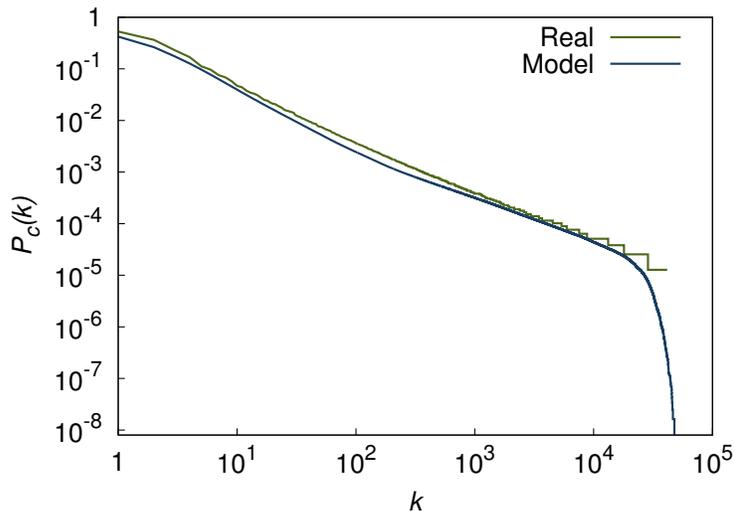}
\caption{\label{omp_pk} Complementary cumulative degree distribution $P_{c} (k)$ of the one-mode projection graph for both the real and the stochastic model networks.}
\end{figure}

In fig.~\ref{omp_pk} we compare Eq.~\eqref{pkonemode} with its stochastic homologous.

\subsection{Weight of a connection and strength of a prime}
The weight of the connection $\omega_{ij}$ between two primes $p_{i}$ and $p_{j}$ is
\begin{equation}\label{weight}
\omega_{ij} = \left\lfloor \frac{N}{p_{i} p_{j}} \right\rfloor.
\end{equation}
This quantity is defined as the number of composites in the bipartite network to which both primes are connected. Such composites must be divisible by both $p_{i}$ and $p_{j}$, i.e. by $p_{i} p_{j}$. Since there are $\left\lfloor N/p_{i} p_{j} \right\rfloor$ such numbers amongst the first $N$ natural numbers, Eq.~\eqref{weight} effectively gives $\omega_{ij}$.

The strength of a prime number is straightforward to obtain from the latter result. By its definition, the only thing to do is adding the weights of the connections to all the other prime numbers in the network, from $p_{1}$ to $p_{\pi \left( \frac{N}{p} \right)}$ (notice that if $p_{i} > N/p$ the weight of the connection is equal to zero). This leads to
\begin{equation}\label{strengthone}
s(p) = \sum \limits_{i=1}^{\pi \left( \frac{N}{p} \right)} \left\lfloor \frac{N}{p p_{i}} \right\rfloor.
\end{equation}
Eq.~\eqref{strengthone} also adds the weight of the self-loop of $p$ if existing (if $p^2 \leq N$).

\subsection{Clustering coefficient}
We have found an expression for the clustering coefficient $ C(p) $ of a prime number inside this graph. This quantity is a real number $ C(p) \in \left[ 0, 1 \right] $ representing the fraction of possible links between the neighbours of $ p $ that actually exist. This coefficient affects many processes in networks such as percolation, dynamic processes, etc. and it is closely related to the small-world property as well as to hidden geometries. In our case, if $ p \geq \sqrt{N} $, it can only be connected to primes $ p_{i} \leq \sqrt{N} $. As the product of two numbers below $ \sqrt{N} $ cannot be greater than $ N $, all the primes $ p_{i} \leq \sqrt{N} $ are connected to each other. Consequently, the clustering coefficient is $ C(p) = 1$ for any $p \geq \sqrt{N} $. However, when $ p \leq \sqrt{N} $, the expression for $ C(p) $ is given by
\begin{equation}\label{clustering}
C(p) = \frac{\left[ \pi \left( p \right) -1 \right] \left[ 2 \left(\pi \left( \frac{N}{p} \right) - 1 \right) - \pi \left( p \right) \right] +
2 \left[ \sum \limits_{j=\pi \left( p \right) + 1}^{\pi(\sqrt{N})} \left( \pi \left( \frac{N}{p_{j}} \right) - j \right) +
 \pi \left( \sqrt{N} \right) - 1 \right]}
{\pi \left( \frac{N}{p} \right) \left[ \pi \left( \frac{N}{p} \right) -1 \right]}.
\end{equation}
To derive Eq.~\eqref{clustering} we need to count the number of connections between the primes to which $p$ is connected. Let us compute several quantities separately.
\begin{itemize}
\item[1.] The number of neighbours of the prime $p$ that we need to consider is not $k(p)$ as given by Eq.~\eqref{om_degree}, but $k \equiv k(p)-1=\pi \left( \frac{N}{p} \right)-1$; since $p \leq \sqrt{N}$, we must correct the fact that Eq.~\eqref{om_degree} is counting the self-loop of prime $p$. The number of possible links amongst these nodes is, allowing the possibility for self-loops to exist,
\begin{equation}\label{Lmax}
L_{max} = \frac{1}{2} k \left( k + 1 \right) = \frac{1}{2} \pi \left( \frac{N}{p} \right) \left[ \pi \left( \frac{N}{p} \right) -1 \right].
\end{equation}

\item[2.] The number of self-loops existing amongst the neighbours of $p$, $L_{sl}$, can be derived easily; a self-loop exists if and only if the corresponding prime is less or equal to $\sqrt{N}$. In addition, $p$ is connected to all such primes, so
\begin{equation}\label{Lsl}
L_{sl} = \pi \left( \sqrt{N} \right) - 1.
\end{equation}
The minus one term corrects the overcount due to the self-loop of prime $p$. This result allows us to simply count the number of links amongst the neighbours of $p$ regardless of self-loops. This calculation is conveniently separated into two more parts.

\item[3.] Links concerning primes less than $p$: let $ p_{i} $ denote any prime less than $p$ (so $i=1,\ldots,\pi \left( p \right) - 1$). Then, if for some prime $p'$ it is true that $ p p' \leq N $, it must be true that $p_{i} p' < N$. In other words, all the $ p_{i} $ are connected to all the primes to which $ p $ is connected. Therefore, we need to count the number of different connections that $ \pi \left( p \right) -1 $ elements can form with $ k $ elements (regardless of self loops, as explained above). We can proceed in the following manner: the first of the $ p_{i} $, $ p_1 $, is connected to $ k-1 $ elements. The second prime, $ p_2 $, forms $ k-2 $ new bonds, since the connection to $p_1$ is not counted again. The elements in this succession can be written as $ k - j $, which allows us to write the corresponding series as
\begin{equation*}
L_{p_{i} < p} = \sum \limits_{j=1}^{\pi \left( p \right) -1} \left( k - j \right)
 = \left( \pi \left( p \right) -1 \right) k - \sum \limits_{j=1}^{\pi \left( p \right) -1} j = \left( \pi \left( p \right) -1 \right) k - \frac{\left(\pi \left( p \right) -1\right)\pi \left( p \right)}{2}.
\end{equation*}
Making now use of the expression for $k$ derived previously yields
\begin{equation}\label{Lless}
L_{p_{i} < p} = \frac{1}{2} \left[ \pi \left( p \right) -1 \right] \left[ 2 \left(\pi \left( \frac{N}{p} \right) - 1 \right) - \pi \left( p \right) \right].
\end{equation}

\item[4.] Links not concerning primes less than $p$: consider any pair of primes $p_{i}$ and $p_{j}$ such that $p_{j} > p_{i} > p$. Then, if $p_{i} p_{j} \leq N$, the chained inequalities $ p p_{i} < p p_{j} < N $ must hold as well. This means that any link between $p_{i}$ and $p_{j}$ (both greater than $p$), is a link amongst neighbours of $p$; in particular, those that we have not counted yet. An easy way to count such links is to count, for every $p_{i}$, the number of $p_{j}$ such that $p_{i} p_{j} \leq N$. For any given $p_{i}$ we see that the value for $p_{j}$ is bounded by $p_{i} < p_{j} \leq N/p_{i}$, so there are $\pi \left( N/p_{i} \right) - \pi \left( p_{i} \right) = \pi \left( N/p_{i} \right) - i$ links to be counted for prime $p_{i}$. The only thing left to do is adding the terms for all the $p_{i}$. Note, however, that the upper bound for $i$ is given by $ i \leq \pi \left( \sqrt{N} \right)$ (if both $p_{i}$ and $p_{j}$ are greater than $\sqrt{N}$, their product cannot belong to the bipartite network). Finally, we can write
\begin{equation}\label{Lgreater}
L_{p_i > p} = \sum \limits_{j=\pi \left( p \right) + 1}^{\pi(\sqrt{N})} \left( \pi \left( \frac{N}{p_{j}} \right) - j \right).
\end{equation}
\end{itemize}
Eq.~\eqref{clustering} is obtained directly by adding Eqs.~(\ref{Lsl}--\ref{Lgreater}) and dividing the result by Eq.~\eqref{Lmax}.
\begin{equation}
C(p) = \frac{L_{sl}+L_{p_i < p}+L_{p_i > p}}{L_{max}}.
\end{equation}

\begin{figure}[h!]
\centering
\includegraphics[scale=0.4]{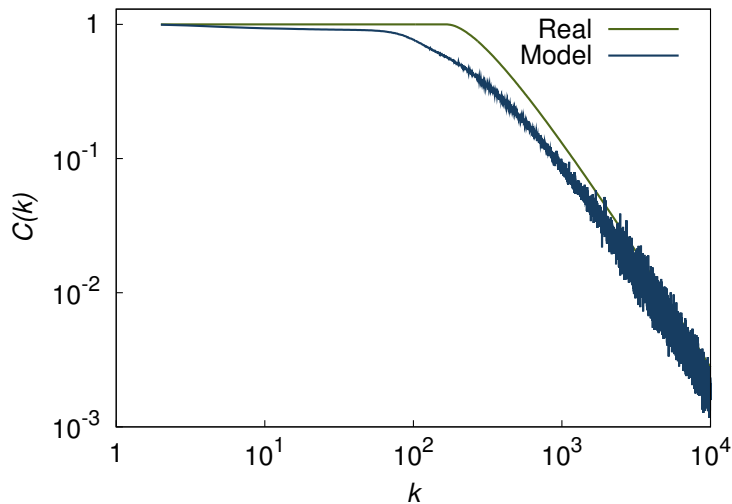}
\caption{\label{omp_ck} Clustering coefficient as a function of the degree $C(k)$.}
\end{figure}

We have obtained a numerical relation between the clustering coefficient $C$ and the degree $k$ as well, which is plotted in Fig.~\ref{omp_ck} with the corresponding measurement on the stochastic model.

\section{Arithmetic functions}
\label{appendixC}
The perspective of number theory that we have presented in this work provides us with a new approach to some arithmetic functions as well. In this section, we present a few results derived from our network representation of natural numbers concerning several of them. We have been able to derive exact and approximated expressions for the prime factors counting function $\omega(n)$ (indeed, we have informally obtained its normal order in accordance with the Hardy-Ramanujan theorem), the sum of the prime divisors of a number $n$ raised to the $r$-th power (which we denote by $\tau_{r}(n)$) and, indirectly, the sum of divisors of $n$ to the $r$-th power, $\sigma_{r}(n)$.

\subsection{Prime factors counting function $\omega(n)$}
In the bipartite network that we have studied, every link connects a prime and a composite. Therefore, counting all the distinct links in the graph (i.e., with no multiplicities) yields the sum of the distinct prime divisors of all the composites up to $N$, 
\begin{equation}
\sum \limits_{n ~ \text{composite} ~ \leq N} \omega(n) = \pi \left( N \right) \sum \limits_{k_{p}} k_{p} P(k_{p}).
\end{equation}
Since $\omega(p)=1$ for any prime, we can extend the latter sum to all $n \in [2,N]$ simply as
\begin{equation}
\sum \limits_{n=2}^{N} \omega(n) = \pi \left( N \right) \left[ 1 + \sum \limits_{k_{p}} k_{p} P(k_{p}) \right].
\end{equation}
Expanding the sum over $k_{p}$ gives
\begin{equation}
\pi \left( N \right) \sum \limits_{k_{p}} k_{p} P(k_{p}) = \pi \left( \frac{N}{2} \right) - \pi \left( \frac{N}{3} \right) + 2 \left[ \pi \left( \frac{N}{3} \right) - \pi \left( \frac{N}{4} \right) \right ] + \cdots = \sum \limits_{k \geq 2} \pi \left( \frac{N}{k} \right).
\end{equation}
We can find an upper limit for the sum in the latter expression considering that $\pi \left( N/k \right) > 0 \Leftrightarrow  N/k \geq 2$, so only the terms with $k \leq \left \lfloor N/2 \right \rfloor$ need to be added. We are finally led to the interesting identity
\begin{equation}\label{omega_sum}
\sum \limits_{n=2}^{N} \omega(n) = \sum \limits_{k = 1}^{\left \lfloor N/2 \right \rfloor} \pi \left( \frac{N}{k} \right).
\end{equation}
The arithmetic function $\omega(n)$ is given in terms of Eq.~\eqref{omega_sum} as the difference between two consecutive sums, i.e. between the sums up to $n$ and $n-1$,
\begin{equation}\label{omega_pi}
\omega(n) = \sum \limits_{k = 1}^{\left \lfloor n/2 \right \rfloor} \left[ \pi \left( \frac{n}{k} \right) - \pi \left( \frac{n-1}{k} \right) \right].
\end{equation}

\subsection{Sum of prime factors of $n$ raised to the $r$-th power $\tau_{r}(n)$}
A further analysis of Eq.~\eqref{omega_pi} reveals that $\phi(k;n) \equiv \pi \left( \frac{n}{k} \right) - \pi \left( \frac{n-1}{k} \right)$ gives
\begin{equation}\label{phi}
\phi(k;n) = \left\{
\begin{array}{ll}
1 & \text{if}~\frac{n}{k}~\text{is prime} \\
0 & \text{otherwise}\\
\end{array}
\right.
\end{equation}
This result allows us to write an expression for $\tau_{r}(n)$, which we define as the sum of prime divisors of $n$ raised to the $r$-th power,
\begin{equation}\label{tau}
\tau_{r}(n) = \sum \limits_{k = 1}^{\left \lfloor n/2 \right \rfloor} \left( \frac{n}{k} \right)^{r} \phi(k;n),
\end{equation}
so $\omega(n) = \tau_{0}(n)$. However, we need to prove Eq.~\eqref{phi}.

First notice that $\phi(k;n)$ is equal to the number of primes in the interval $p \in \left( \frac{n-1}{k}, \frac{n}{k} \right]$. Let us thus count the number of integers in the interval. Suppose that $\frac{n}{k} \notin \mathbb{N}$. Then,
\begin{equation}\label{nk_notinteger}
\frac{n}{k} = \frac{qk + r}{k} > q,
\end{equation}
with $q = \left \lfloor \frac{n}{k} \right \rfloor \in \mathbb{N}$ and $r \geq 1 \in \mathbb{N}$. In addition, we have
\begin{equation}\label{n-1k_notinteger}
\frac{n-1}{k} = \frac{qk + r - 1}{k} \geq q.
\end{equation}
Even though $\frac{n-1}{k}$ can be an integer (if $r=1$), it does not belong to the interval $\left( \frac{n-1}{k}, \frac{n}{k} \right]$, so every number in the interval is greater than $q$. We thus conclude that if $\frac{n}{k} \notin \mathbb{N}$ there are no integers (and therefore no primes) in the interval ($\frac{n}{k} \notin \mathbb{N} \Rightarrow \phi(k;n) = 0$).

On the other hand, if $\frac{n}{k} \in \mathbb{N}$, Eq.~\eqref{nk_notinteger} reads
\begin{equation}\label{nk_integer}
\frac{n}{k} = q,
\end{equation}
while Eq.~\eqref{n-1k_notinteger} becomes
\begin{equation}\label{n-1k_integer}
\frac{n-1}{k} = \frac{\left( q - 1 \right) k + k - 1}{k} \geq q - 1.
\end{equation}
In this case, we see that every number in the interval $x \in \left( \frac{n-1}{k}, \frac{n}{k} \right]$ lies between $q-1<x<q$ (and, hence, they cannot be integers) except for $x=\frac{n}{k} \in \mathbb{N}$. We can thus conclude that, if $\frac{n}{k}$ is prime, $\pi \left( \frac{n}{k} \right) = \pi \left( \frac{n-1}{k} \right) + 1$ and, therefore, $\phi(k;n) = 1$. Notice, however, that even though $\frac{n}{k} \in \mathbb{N}$, if it is not prime, $\pi \left( \frac{n}{k} \right) = \pi \left( \frac{n-1}{k} \right) \Leftrightarrow \phi(k;n) = 0$.

\subsection{Approximation of $\tau_{r}(n)$}
We can derive an approximation of Eq.~\eqref{omega_pi} exchanging the sum for an integral and making use of the prime number theorem (from Eq.~\eqref{pnt}, we see that $\pi \left( x \right) \sim x/\ln x$),
\begin{equation}\label{omega_approx_1}
\begin{aligned}
\omega(n) = \tau_{0} (n) &= \sum \limits_{k = 1}^{\left \lfloor n/2 \right \rfloor} \left[ \pi \left( \frac{n}{k} \right) - \pi \left( \frac{n-1}{k} \right) \right] \sim \int \limits_{1}^{n/2} \left[ \frac{n}{k \ln \frac{n}{k}} - \frac{n-1}{k \ln \frac{n-1}{k}} \right] \text{d}k \\
&\sim \int \limits_{1}^{n/2} \frac{\text{d}k}{k \ln \frac{n}{k}} = \int \limits_{2}^{n} \frac{\text{d}p}{p \ln p} = \ln \ln n - \ln \ln 2.
\end{aligned}
\end{equation}
The latter expression yields, according to the Hardy-Ramanujan theorem, the normal order of $\omega (n)$.

By the same line of reasoning, we can approximate any of the $\tau_{r}(n)$ for $r>0$,
\begin{equation}\label{tau_approx_1}
\tau_{r}(n) \sim \int \limits_{1}^{n/2} \left( \frac{n}{k} \right)^{r} \frac{\text{d}k}{k \ln \frac{n}{k}} = \int \limits_{2}^{n} \frac{p^{r-1}}{\ln p} \text{d}p = \int \limits_{2^{r}}^{n^{r}} \frac{\text{d}q}{q} = \text{li}(n^{r}) - \text{li}(2^{r}).
\end{equation}
Eq.~\eqref{tau_approx_1} yields a very interesting result; in the particular case of $r=1$, we see that $\tau_{1}(n) \sim \text{Li}(n) \sim \pi(n)$, i.e., the sum of the distinct prime factors of $n$ is close to the the number of primes up to $n$.

\subsection{Sum of divisors of $n$ raised to the $r$-th power $\sigma_{r}(n)$}
The proof of Eq.~\eqref{phi} can be used to find an expression for $\sigma_{r}(n)$, defined as the sum of the divisors of $n$ raised to the $r$-th power. Indeed, using Eqs.~(\ref{nk_notinteger},\ref{n-1k_notinteger}) we see that, if $\frac{n}{k} \notin \mathbb{N} \Rightarrow \left \lfloor \frac{n}{k} \right \rfloor = q = \left \lfloor \frac{n-1}{k} \right \rfloor$. On the other hand, if $\frac{n}{k} \in \mathbb{N} \Rightarrow \left \lfloor \frac{n}{k} \right \rfloor = q$ but $ \left \lfloor \frac{n-1}{k} \right \rfloor = q - 1$. If we define $\psi(k;n) \equiv \left \lfloor \frac{n}{k} \right \rfloor - \left \lfloor \frac{n-1}{k} \right \rfloor$ we can write
\begin{equation}\label{psi}
\psi(k;n) = \left\{
\begin{array}{ll}
1 & \text{if}~\frac{n}{k} \in \mathbb{N} \\
0 & \text{otherwise}\\
\end{array}
\right.
\end{equation}
As a consequence, we can easily sum all the divisors of $n$ raised to any power $r$ simply as
\begin{equation}\label{sigma}
\sigma_{r}(n) = \sum \limits_{k=1}^{n} \left( \frac{n}{k} \right)^{r} \psi(k;n) = \sum \limits_{k=1}^{n} k^{r} \psi(k;n).
\end{equation}
The reason why the two sums in Eq.~\eqref{sigma} are equivalent is that, if $\frac{n}{k} \in \mathbb{N}$, both $\frac{n}{k}$ and $k$ divide $n$.

\section{The prime counting function in the stochastic model}
\label{appendixD}
In this section we derive an expression for $P_{N}$, the probability that $N$ is prime in a network chosen at random from the set of all networks of size greater or equal to $N$ generated by our model. Since this probability cannot depend on numbers that join the network after $N$, we only need to study these networks up to $N$ in the calculation. We can describe the state of a particular realization up to $N$ using the set of dichotomous random variables $\left( n_{2},\ldots,n_{N} \right)$, where
\begin{equation}
n_{k} = \left\{
\begin{array}{ll}
1 & \text{if}~k~\text{is prime} \\
0 & \text{otherwise}\\
\end{array}
\right.
\mbox{ with }
k=2,\cdots,N
\end{equation}
This allows us to write $P_{N}$ as
\begin{equation}
P_{N} = \left< n_{N} \right> = \sum \limits_{n_{2}=0}^{1} \cdots \sum \limits_{n_{N}=0}^{1} n_{N} \rho \left( n_{2},\ldots,n_{N} \right),
\end{equation}
where $\langle \cdot \rangle$ denotes the statistical average and $\rho \left( n_{2},\ldots,n_{N} \right)$ is the joint probability of the particular sequence $\left( n_{2},\ldots,n_{N} \right)$. It is convenient to define its characteristic function
\begin{equation}\label{phat}
\hat{\rho} \left( z_{2},\ldots,z_{N} \right) \equiv \sum \limits_{n_{2}=0}^{1} \cdots \sum \limits_{n_{N}=0}^{1} z_{2}^{n_{2}} \ldots z_{N}^{n_{N}} \rho \left( n_{2},\ldots,n_{N} \right).
\end{equation}
$P_{N}$ can be derived from this expression as
\begin{equation}\label{derivative}
P_{N} = \left. \frac{\partial \hat{\rho}}{\partial z_{N}} \right|_{z_{2}=z_3=\cdots=z_N=1}.
\end{equation}
The set of random variables $\left( n_{2},\ldots,n_{N} \right)$ defines a sequence of causal variables, in the sense that $n_i$ only depends on $n_j$ with $j<i$. This implies that $\rho \left( n_{2},\ldots,n_{N} \right)$ satisfies the following Chapman-Kolmogorov equation
\begin{equation}
\rho \left( n_{2},\ldots,n_{N} \right)=\rho \left( n_{2},\ldots,n_{N-1} \right) \mbox{Prob}\{ n_{N} | n_{2},\ldots,n_{N-1}\},
\label{chapman-kolmogorov}
\end{equation}
with $N\ge 3$ and the initial condition $\rho(n_2=1)=1$. The conditional probability that $N$ is prime given the sequence $\left( n_{2},\ldots,n_{N-1} \right)$ is the probability that $N$ does not connect to any of the existing primes below $\sqrt{N}$, that is
\begin{equation}
\mbox{Prob}\{ n_{N} =1| n_{2},\ldots,n_{N-1}\}=\prod \limits_{i=2}^{\left \lfloor \sqrt{N} \right \rfloor} \left( 1 - \frac{1}{i} \right)^{n_{i}},
\end{equation}
and $\mbox{Prob}\{ n_{N} =0| n_{2},\ldots,n_{N-1}\}=1-\mbox{Prob}\{ n_{N} =1| n_{2},\ldots,n_{N-1}\}$. Plugging this expression in Eq.~(\ref{chapman-kolmogorov}) and then to Eq.~(\ref{phat}) leads to the following recurrence relation
\begin{equation}\label{recurrence}
\hat{\rho} \left( z_{2},\ldots,z_{N} \right) = \hat{\rho} \left( z_{2},\ldots,z_{N-1} \right) + \left( z_{N} - 1 \right) \hat{\rho} \left( z_{2} \alpha_{2},\ldots,z_{\left \lfloor \sqrt{N} \right \rfloor} \alpha_{\left \lfloor \sqrt{N} \right \rfloor},z_{\left \lfloor \sqrt{N} \right \rfloor + 1}, \ldots, z_{N-1} \right),
\end{equation}
where we have defined the compact notation $\alpha_{i} \equiv 1 - 1/i$. Finally, by making use of Eq.~\eqref{derivative} we obtain
\begin{equation}
P_{N} = \hat{\rho} \left( \alpha_{2},\ldots,\alpha_{\left \lfloor \sqrt{N} \right \rfloor} \right).
\label{P_N}
\end{equation}
From Eq.~(\ref{recurrence}) it is clear that the random variables $\left( n_{2},\ldots,n_{N-1} \right)$ are not statistically independent. This implies that the exact solution of the problem can only be obtained by solving Eq.~(\ref{recurrence}) and plugging the solution in Eq.~(\ref{P_N}), a task that is, currently, beyond our mathematical skills. Nevertheless, it is possible to derive a very accurate mean field approximation. We start by expanding $\hat{\rho} \left( z_{2},\ldots,z_{\left \lfloor \sqrt{N} \right \rfloor} \right)$ around $z_{1} =z_2=\cdots=z_{\left \lfloor \sqrt{N} \right \rfloor}= 1$ as
\begin{equation}\label{prob_expand}
P_{N} = 1 + \sum \limits_{i=2}^{\left \lfloor \sqrt{N} \right \rfloor} \left.\frac{\partial \hat{\rho}}{\partial z_{i}}\right|_{z_i=1} \beta_{i} + \frac{1}{2!}\sum \limits_{i=2}^{\left \lfloor \sqrt{N} \right \rfloor} \sum \limits_{j=2}^{\left \lfloor \sqrt{N} \right \rfloor} \left.\frac{\partial^{2} \hat{\rho}}{\partial z_{i} \partial z_{j}}\right|_{z_i=z_j=1} \beta_{i} \beta_{j} + \frac{1}{3!} \sum \limits_{i=2}^{\left \lfloor \sqrt{N} \right \rfloor} \sum \limits_{j=2}^{\left \lfloor \sqrt{N} \right \rfloor} \sum \limits_{k=2}^{\left \lfloor \sqrt{N} \right \rfloor} \left.\frac{\partial^{3} \hat{\rho}}{\partial z_{i} \partial z_{j} \partial z_{k}}\right|_{z_i=z_j=z_k=1} \beta_{i} \beta_{j} \beta_{k} + \cdots
\end{equation}
where we have used the convenient notation $\beta_{i} \equiv \alpha_{i} - 1=-1/i$. All terms in the latter expansion that involve a derivative of order higher than one with respect to any of the $z_{i}$ are null, since $n_{i} \left( n_{i} - 1 \right) = 0$ ($n_{i}$ is either 0 or 1). Using this fact and the properties of generating functions, we can rewrite Eq.~\eqref{prob_expand} as
\begin{equation}\label{prob_expand2}
P_{N} = 1 + \sum \limits_{i} \left< n_{i} \right> \beta_{i} + \sum \limits_{i < j} \left< n_{i} n_{j} \right> \beta_{i} \beta_{j} + \sum \limits_{i < j < k} \left< n_{i} n_{j} n_{k} \right> \beta_{i} \beta_{j} \beta_{k} + \cdots+
\langle n_2 n_3 \cdots n_{\left \lfloor \sqrt{N} \right \rfloor}\rangle \beta_2 \beta_3 \cdots \beta_{\left \lfloor \sqrt{N} \right \rfloor}.
\end{equation}
Despite the fact that random variables $n_i$ are not statistically independent, in most of the cases they are conditionally independent. For instance, let us first consider the term $\left< n_{i} n_{j} \right>$ for $i>j$. If $j> \sqrt{i}$ then the only correlation between $n_i$ and $n_j$ is given through their common history, that is, the sequence of primes up to $\sqrt{j}$ and, therefore, they are conditionally independent. In the opposite case, $n_i$ is correlated to $n_j$. However, notice that i) $n_j$ is only one out of $\sqrt{i}$ variables that have a direct influence on $n_i$. ii) The common history between $n_i$ and $n_j$ is even smaller than before and iii) the number of correlated terms for a given $N$ is $\sum_{j=3}^{\sqrt{N}} \sqrt{j} \sim N^{3/4}$ whereas the total number of terms scales as $N^2$. Given these considerations, it is quite reasonable to factorize $\left< n_{i} n_{j} \right> \approx \langle n_i \rangle \langle n_j \rangle=P_i P_j$. A similar analysis can be performed for higher order correlation functions. Under this approximation, Eq.~\eqref{prob_expand2} can be written as
\begin{equation}
P_{N} \approx 1 + \sum \limits_{i} \left< n_{i} \right> \beta_{i} + \sum \limits_{i < j} \left< n_{i} \right> \left< n_{j} \right> \beta_{i} \beta_{j} + \sum \limits_{i < j < k} \left< n_{i} \right> \left< n_{j} \right> \left< n_{k} \right> \beta_{i} \beta_{j} \beta_{k} + \ldots + \left< n_{2} \right> \cdots \left< n_{\left \lfloor \sqrt{N} \right \rfloor} \right> \beta_{2} \cdots \beta_{\left \lfloor \sqrt{N} \right \rfloor}.
\end{equation}
The latter sum can be expressed as
\begin{equation}
P_{N} \approx \sum \limits_{m_{2}=0}^{1} \cdots \sum \limits_{m_{\left \lfloor \sqrt{N} \right \rfloor}=0}^{1} \left( \left< n_{i} \right> \beta_{i} \right)^{m_{i}} = \prod \limits_{i=2}^{\left \lfloor \sqrt{N} \right \rfloor} \left( 1 + \left< n_{i} \right> \beta_{i} \right) = \prod \limits_{i=2}^{\left \lfloor \sqrt{N} \right \rfloor} \left( 1 - \frac{P_{i}}{i} \right).
\end{equation}
Finally, we can write
\begin{equation}
P_{N} \approx e^{\displaystyle{\sum \limits_{i=2}^{\left \lfloor \sqrt{N} \right \rfloor} \ln \left( 1 - \frac{P_{i}}{i} \right)}}.
\end{equation} 
In the limit $N\rightarrow \infty$, the sum in the exponent of the exponential function is dominated by the upper limit and, therefore, it can be approximated as
\begin{equation}
P_{N} \approx e^{\displaystyle{-\sum \limits_{i=2}^{\left \lfloor \sqrt{N} \right \rfloor} \frac{P_{i}}{i} }}.
\label{P_Napprox}
\end{equation}

\section{The Erd\"os-Kac theorem in the stochastic model}
\label{appendixE}
The Erd\"os-Kac theorem states that the quantity $(\omega(N)-\ln{\ln{N}})/\sqrt{\ln{\ln{N}}}$ behaves as a random variable that follows a standard normal distribution. This is known as the fundamental theorem of probabilistic number theory. In our model, this quantity is, indeed, a random variable. In this section, we develop an approximation for the probability that number $N$ in our model has $\omega$ distinct prime factors, $P(\omega | N)$. To do so, we first define the set of dichotomous random variables $(m_2,m_3,\cdots,m_{\left \lfloor \sqrt{N} \right \rfloor})$ as follows
\begin{equation}
m_{k} = \left\{
\begin{array}{ll}
1 & \text{if}~k~\text{is a prime factor of $N$} \\
0 & \text{otherwise}\\
\end{array}
\right.
\mbox{ with }
k=2,\cdots,\left \lfloor \sqrt{N} \right \rfloor
\end{equation}
In terms of these variables, we can write
\begin{equation}
P(\omega| N)=\sum \limits_{n_{2}=0}^{1} \cdots \sum \limits_{n_{\left \lfloor \sqrt{N} \right \rfloor}=0}^{1} \rho(n_2,\cdots,n_{\left \lfloor \sqrt{N} \right \rfloor})
\sum \limits_{m_{2}=0}^{1} \cdots \sum \limits_{m_{\left \lfloor \sqrt{N} \right \rfloor}=0}^{1}
\mbox{Prob}\{ m_2,\cdots,m_{\left \lfloor \sqrt{N} \right \rfloor} |n_2,\cdots,n_{\left \lfloor \sqrt{N} \right \rfloor}\} \delta_{\omega,1+\sum_i m_i},
\end{equation}
where $\delta_{\cdot,\cdot}$ is the Kronecker delta function. The conditional probability of variables $m_i$ satisfies
\begin{equation}
\mbox{Prob}\{ m_2,\cdots,m_{\left \lfloor \sqrt{N} \right \rfloor} |n_2,\cdots,n_{\left \lfloor \sqrt{N} \right \rfloor}\}=
\mbox{Prob}\{ m_2 | n_2\}
\mbox{Prob}\{ m_3 | n_3,m_2\} \mbox{Prob}\{ m_4 | n_4,m_2,m_3\} \cdots,
\end{equation}
with
\begin{equation}
\mbox{Prob}\{ m_j | n_j,m_2,m_3,\cdots,m_{j-1}\}=\delta_{m_j,1} \frac{n_j}{j} \theta\left(\sqrt{\frac{N}{\prod_{i=1}^{j-1} i^{n_i m_i}}} -j\right)+
\delta_{m_j,0} \left[1- \frac{n_j}{j} \theta\left(\sqrt{\frac{N}{\prod_{i=1}^{j-1} i^{n_i m_i}}} -j\right)\right].
\label{P(m)}
\end{equation}
In the latter expression, $\theta(x)$ is the Heaviside step function. Notice that this step function accounts for the fact that $j$ cannot be a prime factor of $N$ if there already exist smaller prime factors such that $j$ is above the square root of the ratio between $N$ and the product of all prime factors smaller than $j$. Dropping this restriction would correspond to evaluate the distribution of a random variable $\hat{\omega}$ that is an upper bound of $\omega$. However, in the limit $N \rightarrow \infty$, since the probability of $j$ being a prime factor decreases as $1/j$, most of the prime factors of $N$ are small numbers for which the argument of the Heaviside function in Eq.~(\ref{P(m)}) is always positive. We then expect that, in such limit, $\hat{\omega} \rightarrow \omega$ and so we can safely drop the Heaviside function in Eq.~(\ref{P(m)}). Under this approximation, the generating function of $P(\omega|N)$ can be written as
\begin{equation}
\hat{P}(z|N)\equiv \sum_{\omega=1}^{\infty} z^{\omega} P(\omega|N)=z \sum \limits_{n_{2}=0}^{1} \cdots \sum \limits_{n_{\left \lfloor \sqrt{N} \right \rfloor}=0}^{1} \rho(n_2,\cdots,n_{\left \lfloor \sqrt{N} \right \rfloor})
\prod_{j=2}^{\left \lfloor \sqrt{N} \right \rfloor} \left[1+(z-1) \frac{n_j}{j} \right],
\end{equation}
and using the same mean field approximation that we used in the previous section, we can write
\begin{equation}
\hat{P}(z|N)=z \prod_{j=2}^{\left \lfloor \sqrt{N} \right \rfloor} \left[1+(z-1) \frac{P_j}{j} \right] \approx z e^{\displaystyle{(z-1)\sum_{j=2}^{\left \lfloor \sqrt{N} \right \rfloor} \frac{P_j}{j}}}.
\end{equation}
We now use Eq.~(\ref{P_Napprox}) to obtain
\begin{equation}
\hat{P}(z|N) \approx z e^{-(z-1) \displaystyle{\ln{P_N}}},
\end{equation}
or, equivalently
\begin{equation}
P(\omega|N)= \frac{P_N}{(\omega-1)!} \left[ -\ln{P_N}\right]^{\omega-1}.
\end{equation}
This is nothing but a Poisson distribution of average $-\ln{P_N}\sim \ln{\ln{N}}$ and standard deviation $\sqrt{ \ln{\ln{N}}}$ which, for large $N$, converges to a Gaussian distribution.

\end{document}